\documentclass[11pt]{amsart}
\usepackage{amsmath,amsthm,amscd,amssymb,amsfonts, amsbsy}
\usepackage{latexsym}
\usepackage{txfonts}

\theoremstyle{plain}

\newtheorem{proposition}{Proposition}[section]

\newtheorem{lemma}[proposition]{Lemma}

\theoremstyle{definition}

\theoremstyle{remark}

\newcommand{\eps}{\varepsilon}
\newcommand{\RR}{{\mathbb{R}}}
\newcommand{\NN}{{\mathbb{N}}}
\newcommand{\ZZ}{{\mathbb{Z}}}

\numberwithin{proposition}{section}
\begin{document}

\title[Correction to ``Hardy and BMO spaces"]{Correction to ``Hardy and BMO spaces associated to divergence form
elliptic operators"}

\author[S. Hofmann]{Steve Hofmann}
\address{Department of Mathematics, 
University of Missouri, Columbia, Missouri 65211, USA}
\email{hofmann@math.missouri.edu}
\thanks{The authors are supported by the National Science Foundation}

\author[S. Mayboroda]{Svitlana Mayboroda}
\address{Department of Mathematics, 
Purdue University, W. Layayette, IN 47907-2067, USA}
\email{svitlana@math.purdue.edu}

\maketitle

We present here a correction to an error in our paper \cite{HM}.  We are grateful to Dachun Yang for bringing the error to our attention.   

In \cite{HM} we develop a theory of $H^1$ (Hardy type) and BMO spaces adapted to a second order, divergence form elliptic (aka accretive) operator $L$ in $\mathbb{R}^n$, with complex, $L^\infty$ coefficents.  
In particular we establish the equivalence of several different $H^1$ type spaces, based on membership
in $L^1$ of various square functions and non-tangential maximal functions adapted to $L$ (if $L$
is the Laplacian, then this theory reduces to that of the classical Hardy and BMO spaces).  Among these (and of central importance) is the ``square function" Hardy 
space $H^1_{S_h}$, defined as the completion of the set
\begin{equation}\{f\in L^2(\mathbb{R}^n): S_h f \in L^1(\mathbb{R}^n)\},\end{equation}
with respect to the norm $$\|f\|_{H^1_{S_h}} : = \|S_hf\|_1,$$
where $$S_hf(x):=\left(\iint_{\{(y,t): |x-y|<t\}}|t^2Le^{-t^2L}f(y)|^2\,\frac
{dydt}{t^{n+1}}\right)^{1/2}.$$  We also show that these spaces may be characterized in terms of a molecular decomposition, and it is to this latter point that we
turn our attention in this note.

We recall now the definition of an $L$-adapted molecule.  For a cube $Q\subset\RR^n$ we denote by $l(Q)$ the sidelength of $Q$ and set

\begin{equation}\label{eq1.7}
S_0(Q):=Q, \quad Q_i=2^iQ,\quad \mbox{ and }\quad
S_i(Q):=2^iQ\setminus 2^{i-1}Q \mbox{ for } i=1,2,...,
\end{equation}

\noindent where $2^iQ$ is cube with the same center as $Q$ and
sidelength $2^il(Q)$.

Let $(p_L,\widetilde{p}_L)$ denote the interior of the range of exponents for which the semigroup
$e^{-tL}$ is bounded on $L^p(\mathbb{R}^n)$. 
Given $p\in(p_L,\widetilde{p}_L)$, $\eps>0$ and $M\in\NN, M > n/4$
(we shall henceforth refer to such $(p,\eps,M)$ as {\it allowable}),
a function $m\in L^p(\RR^n)$,
is called a $(p,\eps,M)$-{\it molecule},
if there exists a cube $Q\subset\RR^n$ such that
\begin{equation}\|m\|_{p,\eps ,M,Q} \,:=\,  \sum_{i= 0}^\infty 2^{i(n -
n/p+\eps)}|Q|^{1-1/p}\sum_{k=0}^M\|(\ell(Q)^2L)^{-k}m\|_{L^p(S_i(Q))}\,\,\leq \, 1.
\label{eq3}\end{equation} 
(Remark:  in \cite{HM}, we used the $\ell^\infty$ norm in $i$, rather than the $\ell^1$ norm,
but it turns out to be slightly more technically convenient to use the latter.  In the end, either
norm will work, as we show that any choice of $\varepsilon >0$ defines an equivalent 
Hardy space).

In \cite{HM}, we defined the $L$-adapted molecular Hardy space $H^1_L$ to be the collection of
all  sums $f=\sum \lambda_j m_j$, with $\{\lambda_j\} \subset \ell^1$, where
each $m_j$ is a $(p,\varepsilon,M)$-molecule;  the norm of an element $f$ in this spaces was 
defined to be $inf \sum|\lambda_j|$
where the infimum runs over all such molecular representations of $f$.   We then claimed
that for $H^1_L$ so-defined,
\begin{equation}H^1_L = \widetilde{H}^1_L\,,\end{equation}
for a certain auxiliary molecular space $\widetilde{H}_L^1$ (defined below), which is in fact the space that we actually work with in the paper (indeed, the arguments in \cite {HM} actually give that the 
various other adapted $H^1$ spaces that we consider, such as $H^1_{S_h}$, are equivalent to
$\widetilde{H}^1_L$).
The salient feature of the space $\widetilde{H}_L^1$ is that it has a dense subset
$\widehat{H}^1_L$ (also defined below), on which the molecular norm is obtained by taking an infimum
only over representations that, in particular, converge in (allowable) $L^p$
(a property that we achieve explicitly by truncating in scale).
Such representations are more natural
when dealing with the action of operators for which there may be no weak-type (1,1) theory
(as is the case for us), but only an $L^p$ theory for some range of $p$.  

In fact, our proof of the claimed equivalence (4) is not correct, given our original definition of the molecular space $H^1_L$, and we suspect now that perhaps the claim may not even be true (although it is true in the classical setting, when $L$ is the Laplacian). We shall
therefore modify the definition of $H^1_L$ as follows, taking
$L^p$ convergence of the molecular representations as a starting point.  Given allowable $(p,\eps,M)$ we say that 
\begin{equation}f = \sum \lambda_i m_i\end{equation} 
is a {\bf p-representation} (or more precisely
a $(p,\eps,M)$-representation) of $f$ if  
$\{\lambda_j\}_{j=0}^\infty\in \ell^1$, the $m_i$ are all $(p,\eps,M)$-molecules, {\bf and
the sum in (5) converges in} $L^p$ (it is this $L^p$ convergence that was missing in the original definition of $H^1_L$
in \cite{HM}). 
We then define the $L$-adapted molecular Hardy space $H^1_L$ to be the completion of the space 
\begin{equation*}\mathbb{H}^1_L \equiv\mathbb{H}^1_{L,p,\eps,M}:= 
 \left\{f:\,\mbox{$f$ has a $(p,\varepsilon,M)$-representation}\right\},
\end{equation*}
with respect to the norm 
\begin{equation*}\label{eq1.11}
\|f\|_{\mathbb{H}^1_{L,p,\eps,M}}=
\inf\Bigl\{\sum_{j=0}^\infty|\lambda_j|:\,\,f=\sum_{j=0}^\infty\lambda_j\,
m_j\, \mbox{ is a
$(p,\eps,M)$-representation} \Bigr\}.
\end{equation*}
With this definition, the remainder of the paper \cite{HM}, including the analogue of Fefferman's $H^1$-$BMO$ duality theorem, as well as the equivalence of the molecular space with the various other $H^1$ type spaces (such as $H^1_{S_h}$ defined above), and the fact that the definition of the molecular space is independent of the choice of allowable
$(p,\eps,M)$, is then correct given some minor adjustments that we shall discuss momentarily.  In particular, the claimed equivalence (4) now holds.

Since the proofs in \cite{HM} (in particular, that of Theorem 4.1) show that the molecular space is independent of the choices of allowable $(p,\eps,M)$,
we shall systematically omit the dependence of the spaces and norms on
$(p,\eps,M)$ except when doing so is likely to cause confusion. In the rest of this note, we restrict our attention to allowable $(p,\eps,M)$. 

We recall now the definition of the auxiliary space $\widehat{H}^1_L$. 
Given allowable $(p,\eps,M)$, and $\delta
> 0,$ we say that $f = \sum \lambda_j m_j$ is a
$\delta$-representation (or more precisely, a $(\delta,p,\eps,M)$-representation) of $f$ if $\{\lambda_j \}_{j=0}^\infty \in
\ell^1$ and each $m_j$ is a $(p,\eps,M)$-molecule adapted to a cube
$Q_j$ of side length {\bf at least} $\delta$. We set
$$H^1_{L,\delta,p,\eps,M}(\RR^n)\equiv H^1_{L,\delta}(\RR^n) := \{ f \in L^1(\RR^n): f \textrm{ has a
$(\delta,p,\eps,M)$-representation}\}.$$ Observe that a $\delta$-representation
is also a $\delta^{'}$-representation for all $\delta^{'} < \delta.$
Thus, $H^1_{L,\delta} \subseteq H^1_{L,\delta^{'}}$  for $0<
\delta^{'} < \delta.$ Set $$ \widehat{H}^1_{L,p,\eps,M}(\RR^n) \equiv \widehat{H}^1_{L}(\RR^n):= 
\cup_{\delta>0} H^1_{L,\delta,p,\eps,M}(\RR^n),$$ and define
\begin{multline*}\|f\|_{\widehat{H}^1_{L,p,\eps,M}(\RR^n)} 
\equiv \|f\|_{\widehat{H}^1_{L}(\RR^n)} 
:=\\ \inf \,\left\{ \sum_{j=0}^\infty |\lambda_j| :
f = \sum_{j=0}^\infty \lambda_j m_j \textrm{ is a
$(\delta,p,\eps,M)$-representation for some } \delta > 0
\right\}.\end{multline*}   Let $\widetilde{H}^1_{L,p,\eps,M} \equiv \widetilde{H}^1_L$ be
the completion of $\widehat{H}^1_{L,p,\eps,M}$ with
respect to this norm.  We note that, since a $(\delta,p,\eps,M)$-representation
clearly converges in $L^p$ (indeed this was the reason we introduced this space
in \cite{HM}), we have trivially that $\widehat{H}^1_{L,p,\eps,M}\subseteq \mathbb{H}^1_{L,p,\eps,M}$.

As we have mentioned, the proofs in \cite{HM} remain essentially unchanged, except for some minor modifications that we now discuss.  To begin, we now have an almost trivial proof of
Lemma 3.3, which states that an $L^p$ bounded linear, or non-negative
sublinear, operator $T$, which maps $(p,\eps,M)$-molecules {\bf uniformly} into $L^1$,
extends to a bounded operator from $H^1_L$ to $L^1.$  Indeed, by density, it is enough to consider
$f \in \mathbb{H}^1_L$, where $f = \sum \lambda_j m_j$ is a $p$-representation such that
$$\|f\|_{{H}^1_L(\RR^n)} \approx  \sum_{j=0}^\infty |\lambda_j| .$$  Since the sum converges in $L^p$, and since $T$ is bounded on $L^p$, we have
that at almost every point,
\begin{equation}|T(f)| \leq \sum_{j=0}^\infty|\lambda_j| \,|T(m_j)|,\end{equation} and more precisely, if T is linear,
\begin{equation}T(f) =\sum_{j=0}^\infty\lambda_j \,T(m_j).\end{equation}
The $L^1$ bound follows.

The other modification concerns Theorem 4.1, or rather Lemma 4.2, of which
the Theorem is an immediate corollary.  The former states that 
$$H^1_L = H^1_{S_h}$$ (cf. (1)). The latter should be restated as follows, although its proof remains, for the most part, the same.

\begin{lemma}  We fix allowable $(p,\eps,M)$, and suppose that 
$f\in \widehat{H}^1_{L,p,\eps,M}\subseteq\mathbb{H}^1_{L,p,\eps,M}$.  
Then there is a constant $C$ depending only on $n,p,\eps,M$ and ellipticity such that
\begin{equation}\tag{i}\|f\|_{\widehat{H}^1_{L,p,\eps,M}(\RR^n)}\leq 
C\|f\|_{H^1_{S_h}(\RR^n)}
\end{equation}
and
\begin{equation}\tag{ii}
\frac{1}{C}\|f\|_{H^1_{S_h}(\RR^n)}\leq\|f\|_{\mathbb{H}^1_{L,p,\eps,M}(\RR^n)}\leq
 \|f\|_{\widehat{H}^1_{L,p,\eps,M}(\RR^n)}.
\end{equation}
In addition, there exists a sequence $\{f_k\}\subset L^2 \cap H^1_{L,p,\eps,M}$ such that
\begin{equation}\tag{iii}
f_k\to f \,\, \text{in} \, \, \widehat{H}^1_{L,p,\eps,M}.
\end{equation}
Finally, the space $\widehat{H}^1_{L,p,\eps,M}$ is densely contained in $H^1_{S_h}$.
\end{lemma}
Thus, Lemma 4.2 says in particular that each allowable $\widehat{H}^1_{L,p,\eps,M}$
is dense in $\widetilde{H}^1_{L,p,\eps,M}$ (by definition), as well as in
$H^1_{L,p,\eps,M}$ (trivially, since any finite linear combination of 
molecules is a $\delta$-representation), and in $H^1_{S_h}$, with equivalence of all of the 
various norms, so Theorem 4.1 follows.

We now sketch the proof of the Lemma, or rather just the modifications that need to be made to
the proof of the original version of  Lemma 4.2 in \cite{HM}.
Estimate (ii) is easy:  the proof of the analogous part of the original Lemma
shows that $S_h$ maps $(p,\eps,M)$-molecules uniformly into $L^1$.  Moreover, $S_h$ is bounded on all allowable $L^p$ \cite{Au}.  Thus, (6) holds with $T = S_h$, and the first inequality in
(ii) follows immediately; the second inequality in (ii) is trivial.

Next, we observe that (iii) is trivial, unless $p<2$ (otherwise, just take $f_k \equiv f,$ 
for every $k$).  If $p<2$, then let $f = \sum_i \lambda_i m_i$ be a $(\delta,p,\eps,M)$-representation of $f$, and for $\eta<\delta$, set \begin{equation}
f_\eta := e^{-\eta^2L} f = \sum_i\lambda_i e^{-\eta^2L}m_i =:\sum_i
\lambda_i m^\eta_i\end{equation} (in the second equality, we have used (7) with $T = e^{-\eta^2L}$). 
By $L^p\to L^2$ hypercontractivity
of the semigroup \cite{Au}, $f_\eta \in L^2$.  Moreover, it is  a routine
matter to verify 
the following: \begin{itemize} \item there is a uniform constant $C_0$ such that for each $i$,
$C_0^{-1} m_i^\eta$ is a $(p,\eps,M)$-molecule adapted to the same cube $Q_i$ as $m_i$.  Thus,
the last sum in (8) is a $(\delta,p,\eps,M)$-representation of $f_\eta;$  \item
$m_i^\eta\to m_i$, as $\eta \to 0$, in the $(p,\eps,M)$-molecular norm defined in (3) (here we write the convergent sum in (3)  as a finite sum plus an arbitrarily small error, and then use that $e^{-\eta^2L}\to I,$
as $\eta \to 0$,
in all allowable $L^p$, to treat each term in the finite sum).  \end{itemize}
Taken together, these two facts readily imply that $f_\eta \to f$, as $\eta \to 0$,
in $\widehat{H}^1_{L,p,\eps,M}$ and thus also in $H^1_{S_h}$ by (ii). 
Consequently, $f\in H^1_{S_h},$ by definition, so that
$$\widehat{H}^1_{L,p,\eps,M}\subset H^1_{S_h}.$$
We now show that this containment is dense.  By density of $L^2 \cap H^1_{S_h}$ in $H^1_{S_h}$,
it will be enough to show that for every $f \in L^2 \cap H^1_{S_h}$, there is a sequence $\{f_k\}\in 
\widehat{H}^1_{L,p,\eps,M}$ such that $f_k \to f$ in $H^1_{S_h}.$  For allowable $p_1<p_2$,
we have that a $(p_2,\eps,M)$-molecule is also a $(p_1,\eps,M)$-molecule, therefore
$$\widehat{H}^1_{L,p_2,\eps,M} \subset \widehat{H}^1_{L,p_1,\eps,M},\quad p_1<p_2.$$
Thus, it is enough to establish the claimed density for $p\geq 2$.   In the original proof of Lemma 4.2 in \cite{HM}, we showed that there exist $f_N \to f $ in $L^2$, with $f_N \in \widehat{H}^1_{L,p,\eps,M}$
for {\it all} allowable $p$, and moreover that $\{f_N\}$ is a Cauchy sequence in
$\widehat{H}^1_{L,p,\eps,M}.$  Thus, $S_h(f_N-f) \to 0$ in $L^2$ and
by Lemma 4.2 (ii) above $$\lim_{N,N'\to\infty}\|S_h(f_N-f_{N'})\|_1 =   0.$$
Taking subsequences, we have that $S_h(f_{N_k}-f) \to 0$ a.e., and
$$\lim_{k' \to \infty} S_h(f_{N_k} - f_{N_{k'}}) = S_h (f_{N_k} - f),\quad \text{a.e.}.$$
Consequently, for any given $\eta > 0$,
$$\|S(f_{N_k} - f)\|_1 =\int \liminf_{k'\to\infty} S_h(f_{N_k} - f_{N_{k'}})\leq 
\liminf_{k'\to\infty}\int S_h(f_{N_k} - f_{N_{k'}})< \eta,$$
for $k$ chosen large enough, and the alleged density follows.

It remains only to discuss the modifications required to prove (i).  Let
$f = \sum \lambda_i m_i$ be a $(\delta,p,\eps,M)$-representation.
For allowable $p$, by $L^p$ functional calculus we have that
$$f = C_M \int_0^\infty \left( t^2 L e^{-t^2L} \right)^{M+2} \!f\, \,\frac{dt}{t} = \lim_{N\to \infty} 
\int_{1/N}^N \left( t^2 L e^{-t^2L} \right)^{M+2}\! f\, \,\frac{dt}{t},$$
where the limit exists in $L^p$.

As in \cite{HM}, we follow the tent space approach of \cite{CMS}.  We 
define the family of sets ${O}_k:= \{x\in\RR^n:\,S_hf(x)>2^k\}$,
$k\in\ZZ$, and consider ${O}_k^*:=\{x\in\RR^n:\, {\mathcal
M}(\chi_{O_k})>1-\gamma\}$ for some fixed $0<\gamma<1$. Then
$O_k\subset O_k^*$ and $|O_k^*|\leq C(\gamma)|O_k|$ for every
$k\in\ZZ$. Next let $\{Q_k^j\}_j$ be a Whitney decomposition of
$O_k^*$ and $\widehat O_k^*$ be a tent region, that is

\begin{equation*}
\widehat O_k^*:=\{(x,t)\in\RR^n\times(0, \infty):\,{\rm
dist}(x,\,^cO_k^*)\geq t\}.
\end{equation*}

\noindent For every $k\in\ZZ, j\in\NN$ we define

\begin{equation*}T_k^j:=\left(Q_k^j\times(0,\infty)\right)\cap \widehat O_k^* \cap
\,^c\widehat O_{k+1}^*,
\end{equation*}

The sets $T_k^j$ are non-overlapping and cover the half-space, so that
\begin{multline*}f = C_M \int_0^\infty \left( t^2 L e^{-t^2L} \right)^{M+1}\left(\sum_{j\in \NN, k\in \ZZ}\chi_{T_k^j}(\cdot,t)\,\,t^2 L e^{-t^2L} f \right)\frac{dt}{t}\frac{dt}{t}\\=
C_M \int_{1/N}^\infty \left( t^2 L e^{-t^2L} \right)^{M+1}\left(\sum_{j\in \NN, k\in \ZZ}
\chi_{T_k^j}\,\,t^2 L e^{-t^2L} f \right)\frac{dt}{t}\\\quad+\,\, C_M \int_0^{1/N} \left( t^2 L e^{-t^2L} \right)^{M+1}\left(\sum_{j\in \NN, k\in \ZZ}\chi_{T_k^j}\,\,t^2 L e^{-t^2L} f \right)\\=: \widetilde{f}_N + T_N(f).\end{multline*}
Then $T_N(f) \to 0$ in $\widehat{H}^1_{L,p,\eps,M}$, as $N \to \infty$ (this is straightforward, and similar to the proof that $f_\eta \to f$ in $\widehat{H}^1_{L,p,\eps,M}$) (cf. (8)).
Therefore, it is enough to show that $$\sup_N\|\widetilde{f}_N\|_{\widehat{H}^1_{L,p,\eps,M}} \leq C\|f\|_{H^1_{S_h}}.$$  In \cite{HM}, this is done in the original proof of Lemma 4.2 for $f_N$, which is defined in the same way as $\widetilde{f}_N$, but with a doubly truncated integral ($\int_{1/N}^N...dt/t$).  
The proof of this fact for $\widetilde{f}_N$ is exactly the same as that for $f_N$, once we show that one can interchange the order of summation and integration, even without truncating the $t$ integral at infinity.
That is, it suffices to prove that \begin{equation}
\widetilde{f}_N = C_M \sum_{j\in \NN, k\in \ZZ}\int_{1/N}^\infty \left( t^2 L e^{-t^2L} \right)^{M+1}\left(
\chi_{T_k^j}\,\,t^2 L e^{-t^2L} f \right)\frac{dt}{t}.
\end{equation}
To this end, recalling that $f \in \widehat{H}^1_{L,p,\eps,M}$ (thus, in particular, in $L^r,\, 1\leq r \leq p$), we choose two allowable exponents $p_1$ and $p_2$, such that $$p_L < p_1 < \min(p,2)\leq \max(p,2) < p_2 < \widetilde{p}_L,$$ and set
$\gamma := n\left(\frac{1}{2} - \frac{1}{p_2}\right)$ and 
$\gamma' := n\left(\frac{1}{p_1} - \frac{1}{2}\right).$   We then have, by hypercontractivity \cite{Au} 
and the fact that the sets $T_k^j$ are non-overlapping,
\begin{multline*}\|\widetilde{f}_N\|_{p_2} \leq \int_{1/N}^\infty \left\| \left( t^2 L e^{-t^2L} \right)^{M+1}\left(\sum_{j\in \NN, k\in \ZZ}
\chi_{T_k^j}\,\,t^2 L e^{-t^2L} f \right)\right\|_{p_2} \frac{dt}{t}\\\leq 
\int_{1/N}^\infty t^{-\gamma}\left\| \left(\sum_{j\in \NN, k\in \ZZ}
\chi_{T_k^j}\,\,t^2 L e^{-t^2L} f \right)\right\|_{2} \frac{dt}{t}\\
\leq C_{N,\gamma} \left(\int_{1/N}^\infty t^{-\gamma}\int_{\RR^n}\left| \left(
\,\,t^2 L e^{-t^2L} f \right)\right|^2 dx\frac{dt}{t}\right)^{1/2}\\
\leq  C_{N,\gamma} \left(\int_{1/N}^\infty t^{-\gamma- 2\gamma'}
\frac{dt}{t}\right)^{1/2}\|f\|_{p_1} = C \|f\|_{p_1}.
\end{multline*}
By dominated convergence, this last argument also shows that
$$\lim_{K\to 0}\left\|\int_{1/N}^\infty \left( t^2 L e^{-t^2L} \right)^{M+1}\left(\sum_{j+|k|> K}
\chi_{T_k^j}\,\,t^2 L e^{-t^2L} f \right)\frac{dt}{t}\right\|_{p_2} = 0.$$
Consequently, to obtain the identity (9), it is enough to prove that, moreover,
\begin{equation}\lim_{K\to 0}\|E_K\|_{p_2} :=\lim_{K\to 0}\left\| \sum_{j +| k| > K}\int_{1/N}^\infty \left( t^2 L e^{-t^2L} \right)^{M+1}\left(
\chi_{T_k^j}\,\,t^2 L e^{-t^2L} f \right)\frac{dt}{t} \right\|_{p_2}  = 0.\end{equation}
As in \cite{HM}, we follow \cite{CMS}
to write
\begin{equation}E_K = \sum_{j +| k| > K}  \lambda_k^j \,\,\widetilde{m}_k^j(N),\end{equation}
where $\lambda_k^j=C_M2^k|Q_k^j|$ (so that $\sum \lambda_k^j \leq C \|S_h f\|_1$) and
\begin{equation*}
\widetilde{m}_k^j(N)=\frac{1}{\lambda_k^j}\int_{1/N}^\infty
(t^2Le^{-t^2L})^{M+1}\left(\chi_{T_k^j}\,t^2Le^{-t^2L}\right)f\,\frac{dt}{t}.
\end{equation*}
Up to a harmless normalization, the $\widetilde{m}_k^j(N)$ are 
$(\widetilde{p},\widetilde{\eps},M)$-molecules, for {\bf all} allowable $\widetilde{\eps}$ and $\widetilde{p}$.   This fact is established in the original
proof of Lemma 4.2 in \cite{HM} for $m_k^j$, which are defined in exactly the same way as
$\widetilde{m}_k^j(N)$, but with a doubly truncated integral ($\int_{1/N}^N...dt/t$).  The proof for the singly truncated integrals considered here is identical.   In addition,  
by definition of $T_k^j$, $\widetilde{m}_k^j(N) = 0$ if $C\ell(Q_k^j) < 1/N,$
so that (11) is a $\widetilde{\delta}$-representation with $\widetilde{\delta} \approx 1/N$.
Consequently, the sum defining $E_K$ converges in all allowable $L^p$, in particular in $L^{p_2}$, so that (10) holds.

To conclude, we note in retrospect that it is not necessary to work explicitly with the
$\delta$-representation space $\widehat{H}^1_{L,p,\eps,M}$.  We could just as well
have proved a version of Lemma 4.2 using instead $\mathbb{H}^1_{L,p,\eps,M}$. 
However, to do so would have required more than the minor revision of the arguments
of \cite{HM} that we have described here.

\end{document}